\documentclass[a4paper,12pt]{amsart}
\usepackage{amssymb}
\usepackage{ifthen}
\usepackage{mathrsfs}
\nonstopmode
\numberwithin{equation}{section}
\newenvironment{thm}[1][]{%
\refstepcounter{equation}%
\medskip%
\noindent%
{\bf Theorem \arabic{section}.\arabic{equation}}%
\ifthenelse{\equal{#1}{}}{}{ (#1)}%
{\bf .}
\itshape}{\medskip}%{\vskip 8pt}

\newenvironment{cor}[1][]{%
\refstepcounter{equation}%
\medskip%
\noindent%
{\bf Corollary \arabic{section}.\arabic{equation}}%
\ifthenelse{\equal{#1}{}}{}{ (#1)}%
{\bf .}
\itshape}{\medskip}%{\vskip 8pt}

\newenvironment{lem}[1][]{%
\refstepcounter{equation}%
\medskip%
\noindent%
{\bf Lemma \arabic{section}.\arabic{equation}}%
\ifthenelse{\equal{#1}{}}{}{ (#1)}%
{\bf .}
\itshape}{\medskip}%{\vskip 8pt}

\newenvironment{prop}[1][]{%
\refstepcounter{equation}%
\medskip%
\noindent%
{\bf Proposition \arabic{section}.\arabic{equation}}%
\ifthenelse{\equal{#1}{}}{}{ (#1)}%
{\bf .}
\itshape}{\medskip}%{\vskip 8pt}

\newenvironment{example}[1][]{%
\refstepcounter{equation}%
\medskip%
\noindent%
{\bf Example \arabic{section}.\arabic{equation}}%
\ifthenelse{\equal{#1}{}}{}{ (#1)}%
{\bf .}
}{\medskip}%{\vskip 8pt}

\newenvironment{rem}[1][]{%
\refstepcounter{equation}%
\medskip%
\noindent%
{\bf Remark \arabic{section}.\arabic{equation}}%
\ifthenelse{\equal{#1}{}}{}{ (#1)}%
{\bf .}
}{\medskip}%{\vskip 8pt}

\newenvironment{pf}[1][]{%
 \vskip 3mm
 \noindent
 \ifthenelse{\equal{#1}{}}%
  {{\slshape Proof. }}%
  {{\slshape #1.} }%
 }%
{\qed\bigskip}

\newcounter{alphabet}
\newcounter{tmp}
\newenvironment{Thm}[1][]{\refstepcounter{alphabet}%
\bigskip%
\noindent%
{\bf Theorem \Alph{alphabet}}%
\ifthenelse{\equal{#1}{}}{}{ (#1)}%
{\bf .}
\itshape}{\vskip 8pt}

\newcounter{minutes}
\divide\time by 60
\newcounter{hours}
\multiply\time by 60
\addtocounter{minutes}{-\time}
\begin{document}
\bibliographystyle{amsplain}
\title{
On $\mu$-conformal homeomorphisms and boundary correspondence
}

\author[V.~Gutlyanski\u\i]{Vladimir Gutlyanski\u\i}
\address{Institute of Applied Mathematics and Mechanics, 
NAS of Ukraine, ul. Roze Luxemburg 74, 83114, Donetsk, Ukraine}
\email{gutlyanskii@iamm.ac.donetsk.ua}
\author[K.~Sakan]{Ken-ichi Sakan}
\address{Department of Mathematics, Graduate School of Science, Osaka City University, Sugimoto, Sumiyoshi-ku, Osaka, 558, Japan}
\email{ksakan@sci.osaka-cu.ac.jp}
\author[T.~Sugawa]{Toshiyuki Sugawa}
\address{Graduate School of Information Sciences,
Tohoku University, Aoba-ku, Sendai 980-8579, Japan}
\email{sugawa@math.is.tohoku.ac.jp}
\date{}
\begin{abstract}
We study the boundary correspondence under $\mu$-homeomorphisms
$f$ of the open upper half-plane onto itself.
Sufficient conditions are given for $f$ to
admit a homeomorphic extension to the closed half-plane with
prescribed boundary regularity.
The proofs are based on the modulus estimates for semiannuli
in terms of directional dilatations of $f$ which might be
of independent interest.
\end{abstract}

\maketitle
%\tableofcontents

\newcommand{\co}{{\overline{\operatorname{co}}}}
\newcommand{\A}{{\mathscr A}}
\newcommand{\B}{{\mathscr B}}
\newcommand{\F}{{\mathscr F}}
\newcommand{\G}{{\mathcal G}}
\newcommand{\Perron}{{\mathscr P}}
\newcommand{\K}{{\mathscr K}}
\newcommand{\es}{{\mathscr S}}
\newcommand{\spl}{{\mathscr S_{\operatorname p}}}
\newcommand{\U}{{\mathscr U}}
\newcommand{\LU}{{\mathscr{LU}}}
\newcommand{\ZF}{{\mathscr{ZF}}}
\newcommand{\R}{{\mathbb R}}
\newcommand{\C}{{\mathbb C}}
\newcommand{\X}{{\mathscr X}}
\newcommand{\PP}{{\mathcal P}}
\newcommand{\h}{{\mathbb H}}
\newcommand{\uhp}{{\mathbb H}}
\newcommand{\Z}{{\mathbb Z}}
\newcommand{\N}{{\mathscr N}}
\newcommand{\M}{{\mathscr M}}
\newcommand{\SCC}{{\mathscr{SCC}}}
\newcommand{\CC}{{\mathscr C}}
\newcommand{\st}{{\mathscr{SS}}}
\newcommand{\D}{{\mathbb D}}
\newcommand{\sphere}{{\widehat{\mathbb C}}}
\renewcommand{\Im}{{\operatorname{Im}\,}}
\newcommand{\Aut}{{\operatorname{Aut}}}
\renewcommand{\Re}{{\operatorname{Re}\,}}
\newcommand{\kernel}{{\operatorname{Ker}}}
\newcommand{\ord}{{\operatorname{ord}}}
\newcommand{\id}{{\operatorname{id}}}
\newcommand{\mob}{{\text{\rm M\"{o}b}}}
\newcommand{\Int}{{\operatorname{Int}\,}}
\newcommand{\Sign}{{\operatorname{Sign}}}
\newcommand{\diam}{{\operatorname{diam}\,}}
\newcommand{\inv}{^{-1}}
\newcommand{\area}{{\operatorname{Area}}}
\newcommand{\iu}{{\operatorname{i}}}
\newcommand{\eit}{{e^{i\theta}}}
\newcommand{\dist}{{\operatorname{dist}}}
\newcommand{\arctanh}{{\operatorname{arctanh}\,}}
\newcommand{\remark}{\vskip .3cm \noindent {\sl Remark.} \@}
\newcommand{\remarks}{\vskip .3cm \noindent {\sl Remarks.} \@}
\newcommand{\ucv}{{\operatorname{UCV}}}
\newcommand{\PSU}{{\operatorname{PSU}}}
\newcommand{\PSL}{{\operatorname{PSL}}}
\newcommand{\belt}{{\operatorname{Belt}}}
\newcommand{\qc}{{\operatorname{QC}}}
\newcommand{\qs}{{\operatorname{QS}}}
\newcommand{\Te}{{\stackrel{\mathrm{T}}{\sim}}}
\newcommand{\teich}{{\operatorname{Teich}}}
\newcommand{\T}{{\operatorname{T}}}
\newcommand{\Th}{{\hat T}}
\newcommand{\Sh}{{\hat S}}
\newcommand{\IR}{{\operatorname{I}}}
\newcommand{\OR}{{\operatorname{O}}}
\newcommand{\Be}{{\mathrm B}}
\newcommand{\sdiam}{{\operatorname{diam}^\sharp}}% spherical diameter
\newcommand{\sd}{{d^\sharp}}% spherical distance
\renewcommand{\mod}{{\,\operatorname{mod}\,}}
\newcommand{\plm}{{m}}% plane area element (2-dim. Lebesgue measure)
\newcommand{\loc}{{\rm{loc}}}% local
\newcommand{\aand}{{\quad\text{and}\quad}}
\newcommand{\db}{{\bar{\partial}}}% d-bar
\newcommand{\dz}{{\bar{z}}}% z-bar
\newcommand{\be}{\begin{equation}}
\newcommand{\ee}{\end{equation}}
\newcommand{\esssup}{\mathop{\operatorname{ess\,sup}\,}\limits}

\section{Introduction and main theorems}

%The regularity problems have been extensively
%studied and they have found many applications, see, e.g.,
%\cite{BJ94}, \cite{BJ02}, \cite{Car61}, \cite{RW65}.
Let $f$ be a quasiconformal homeomorphism of
the open upper half-plane $\h$ onto itself.
A.~Beurling and L.~Ahlfors \cite{BA56} have shown that $f$ admits
a homeomorphic extension to the closed half-plane and that
its boundary map can be characterized as
a quasisymmetric function on $\R$ if it fixes the point at infinity.
Moreover, they remark that the boundary map
need not be locally absolutely continuous on $\R.$
In \cite{Car61}, L.~Carleson proved that the condition
$$%\be\label{eq:carl}
\int_0^1\frac{\eta(s)}{s}dt<\infty,
\quad
\eta(s)=\esssup_{0<\Im z<s}|\mu_f(z)|,~\mu_f=\frac{f_\dz}{f_z},
$$
implies that the boundary correspondence $t\mapsto f(t)$ has a continuous
derivative $f'(t)$ for a quasiconformal self-mapping $f$ of $\uhp$
with $f(\infty)=\infty.$

M.~Brakalova and J.~Jenkins \cite{BJ02}
have extended Carleson's theorem to a class
of $\mu$-conformal homeomorphisms.
Here and hereafter, a {\it Beltrami coefficient} $\mu$ on a domain $\Omega$
will mean a complex-valued measurable function such that
$|\mu|<1$ a.e.~in $\Omega$
and a homeomorphism $f:\Omega\to\Omega'$ will be called $\mu$-conformal
if $f$ belongs to the Sobolev space $W_\loc^{1,1}(\Omega)$
and satisfies the Beltrami equation
$$%\be\label{eq:belt}
f_\dz=\mu f_z \quad\text{a.e.~in~}\Omega.
$$
When $\|\mu\|_\infty<1$ the function $f$ is called {\it quasiconformal}.

For the present and later use, we set
$$
A(z_0;r,R)=\{z\in\C; r\le|z-z_0|\le R\}
$$
for $z_0\in\sphere=\C\cup\{\infty\}$ and
$0\le r<R\le+\infty.$
Here, $A(\infty; r,R)$ is defined to be $A(0; 1/R, 1/r).$

Brakalova and Jenkins \cite{BJ02} proved the following.

\begin{Thm}[Brakalova-Jenkins \cite{BJ02}]\label{Thm:BJ}
Let $f$ be a $\mu$-conformal self-homeomorphism of the upper half-plane
$\uhp$ for a Beltrami coefficient $\mu.$
Suppose that $f(z)\to\infty$ if and only if $z\to\infty$ in $\uhp$
and that
$$
\iint_{A(t,r,R)\cap\uhp}
\frac{|\mu(z)|^2+|\Re \frac{\bar z-t}{z-t}\mu(z)|}{1-|\mu(z)|^2}
\frac{dxdy}{|z-t|^2},\quad
z=x+iy,
$$
converges as $r\to0+$ for every $t\in\R$ and some $R=R(t)>0.$
Then $f$ extends to a homeomorphism of $\overline\uhp$ in such a
way that the boundary function $f:\R\to\R$ is differentiable everywhere and
$$
f'(t)=\lim_{z\to a\text{~in~}\uhp}\frac{f(z)-f(t)}{z-t}>0,
\quad t\in\R.
$$
Moreover if the above convergence is locally uniform for $t\in\R,$
then $f'$ is continuous.
\end{Thm}

In their theorem, continuity of the function at $\infty$ is assumed.
As we will see later (Example \ref{ex:1}), 
this does not follow from the other assumptions.
We slightly refine their result and state it as a local version.

\begin{thm}\label{thm:1}
Let $\mu$ be a Beltrami coefficient on $\uhp$ which satisfies
the following two conditions for $t\in I,$
where $I$ is an open interval in $\R:$
\begin{enumerate}
%\item
%$\displaystyle\lim_{R\to+\infty}\frac{P_\mu(r_0,R)}{\log R}=0,$
%\medskip
\item
$\displaystyle%
%\mathop{\iint}\limits_{|z-t|<R_0 \atop \Im z>0}
\iint_{A(t; r_1,r_2)\cap\uhp}
\frac{|\mu(z)|^2}{1-|\mu(z)|^2}\frac{dxdy}{|z-t|^2}\to 0$
as $r_1, r_2\to0,$ and
\item
$\Re\displaystyle\iint_{A(t; r_1,r_2)\cap\uhp}\frac{\mu(z)}{(z-t)^2}
\frac{dxdy}{1-|\mu(z)|^2}\to0$
as $r_1, r_2\to0.$
\end{enumerate}
Suppose that there exists a $\mu$-conformal self-homeomorphism $f$
of $\uhp.$
Then it extends to a homeomorphism of $\uhp\cup I$ into $\overline\uhp.$
Furthermore, if $f(I)\subset\R,$
then the boundary function $f:I\to\R$ is differentiable and
$$
f'(t)=\lim_{z\to t\text{~in~}\uhp}\left|\frac{f(z)-f(t)}{z-t}\right|>0
$$
holds for each $t\in I.$
Moreover, if the convergence in $(1)$ and $(2)$ is locally uniform
for $t\in I,$ then the derivative $f'$ is continuous on $I.$
\end{thm}

A proof of this theorem will be given in the final section
as well as those of the other results in the present section.
Though the proof of Theorem \ref{thm:1}
can be done in the same way as in \cite{BJ02},
we will supply a detailed account for the continuous extension of
the mapping $f$ since it is largely omitted in \cite{BJ02}.
To this end, in Section 3,
we give an explicit estimate for the modulus of continuity
of the boundary mapping at a given point.
As preparations, we introduce the notion of semiannulus and give
concrete estimates of the modulus of a semiannulus and related quantities
in Section 2.
We believe that these estimates and methods will be useful
in various other problems.
Indeed, as applications of the estimates of modulus of continuity,
we propose a couple of related results in the rest of the present section.
 
We recall that a function $f$ is locally Lipschitz continuous
on a subset $X$ of $\C$
if for every compact subset $E$ of $X,$ there exists a constant $C=C(E)$
such that the inequality $|f(z)-f(z_0)|\le C|z-z_0|$ holds for 
$z_0, z\in E.$
The following result gives us 
a sufficient condition for a $\mu$-conformal homeomorphism $f:\h\to\h$ to
admit an extension to $\overline{\uhp}$
whose boundary correspondence is
locally Lipschitz continuous.
This can be regarded as a boundary version of Theorem 3.13 in \cite{GS01}.

\begin{thm}\label{thm:2}
Let $f:\h\to\h$ be a $\mu$-conformal homeomorphism for
a Beltrami coefficient $\mu$ on $\uhp$
and let $I$ be an open interval in $\R.$
If there are positive constants $R$ and $M$ such that
\be\label{eq:t2}
\iint_{A(t; r, R)\cap\uhp}
%\mathop{\iint}\limits_{r<|z-t|<R \atop \Im z>0}
\frac{|\mu(z)|^2-\Re \frac{\bar z-t}{z-t}\mu(z)}{1-|\mu(z)|^2}
\frac{dxdy}{|z-t|^2}\leq M
\ee
for every $t\in I$ and $r\in(0,R),$
then  $f$ is extended
to a homeomorphism of $\uhp\cup I$ into the closed upper half-plane.
If moreover $f(I)\subset\R,$
the boundary function $f:I\to\R$ is locally Lipschitz continuous on $I.$
\end{thm}

We will say that a function $f$ is 
{\it locally weak H\"older continuous} with exponent $\alpha>0$
in a subset $X$ of $\C$ if for every compact subset $E$ of $X$
and every $0<\alpha'<\alpha$ 
there is a constant $C=C(\alpha',E)>0$ such that
$$
|f(z)-f(z_0)|\leq C|z-z_0|^{\alpha'}
$$
holds whenever $z_0,z\in E.$
In particular, $f$ is called {\it locally weak Lipschitz continuous}
when $\alpha=1.$
We have now the following result, a prototype of which is Theorem 4.5
in \cite{GS01}.

\begin{thm}\label{thm:3}
Let $f:\h\to\h$ be a $\mu$-conformal homeomorphism for a Beltrami coefficient
$\mu.$
For a given open interval $I$ in $\R$ and a number $0<\alpha\le 1,$
we suppose that
\be\label{eq:t3}
\limsup_{r\to 0+}\frac{4}{\pi r^2}
\iint_{A(t; 0, r)\cap\uhp}
%\mathop{\iint}\limits_{|z-t|<R \atop \Im z>0}
\frac{|\mu(z)|^2-\Re \frac{\bar z-t}{z-t}\mu(z)}{1-|\mu(z)|^2}dxdy
\leq \frac{1}{\alpha}-1,
\ee
where the convergence is locally uniform for $t\in I.$
Then $f$ is extended to a homeomorphism of $\uhp\cup I$ into
the closed upper half-plane $\overline\uhp$
and, if $f(I)\subset \R$ in addition,
the boundary correspondence $t\mapsto f(t)$ is 
locally weak H\"older continuous with exponent $\alpha$ on $I.$
\end{thm}

In order to obtain a result of the same type as Theorem A,
we have to impose another condition on $\mu.$

\begin{lem}\label{lem:infty}
Let $f$ be a $\mu$-conformal self-homeomorphism of $\uhp$
for a Beltrami coefficient $\mu$ on $\uhp.$
If
$$
\lim_{R\to+\infty}\frac1{(\log R)^2}
\iint_{A(0; r_0,R)\cap\uhp}
\frac{|\mu(z)|^2-\Re \frac{\bar z}{z}\,\mu(z)}{1-|\mu(z)|^2}
\frac{dxdy}{|z|^2}=0
$$
for some $r_0>0,$ then $f$ extends continuously to the point at infinity.
\end{lem}

We will show the lemma in the last section.
The extended map in the lemma does not necessarily satisfy $f(\infty)=\infty.$
However, $g=M\circ f$ becomes a $\mu$-conformal homeomorphism of $\uhp$
with $g(\infty)=\infty$ for a M\"obius transformation $M$
such that $M(\uhp)=\uhp$ and $M(f(\infty))=\infty.$

We may add the condition in this lemma as well as $I=\R$
to the assumptions in the above theorems to guarantee a homeomorphic
extension of $f$ to the closed upper half-plane.

We also have results analogous to the above theorems 
for the case of the unit disk.
As a sample, we give a variant of Theorem \ref{thm:1}.
In the following, let $T(\zeta; r, R)=\{z\in\D: r<|\frac{z-\zeta}{z+\zeta}|<R\}$
for $\zeta\in\partial\D$ and $0<r<R<+\infty.$

\begin{thm}\label{thm:1d}
Let $\mu$ be a Beltrami coefficient on the unit disk $\D.$
Assume the following two conditions for $\zeta\in\partial\D:$
\begin{enumerate}
\item[(i)]
$\displaystyle\iint_{T(\zeta; r_1,r_2)}
\frac{|\mu(z)|^2}{1-|\mu(z)|^2}\frac{dxdy}{|z^2-\zeta^2|^2}\to 0$
as $r_1, r_2\to0,$ and
\item[(ii)]
$\Re\displaystyle\iint_{T(\zeta; r_1,r_2)}\frac{\zeta^2\mu(z)}{(z^2-\zeta^2)^2}
\frac{dxdy}{1-|\mu(z)|^2}\to0$
as $r_1, r_2\to0.$
\end{enumerate}
Suppose that there exists a $\mu$-conformal self-homeomorphism $f$
of $\D.$
Then it extends to a self-homeomorphism of $\overline\D$
and the boundary map $f:\partial\D\to\partial\D$ is differentiable and
$$
\left|\frac{d}{d\theta}f(e^{i\theta})\right|
=\lim_{z\to e^{i\theta}\text{~in~}\D}
\left|\frac{f(z)-f(e^{i\theta})}{z-e^{i\theta}}\right|>0
$$
holds for each $\theta\in\R.$
Moreover, if the convergence in $(\mathrm{i})$ and $(\mathrm{ii})$
is uniform for $\zeta\in\partial\D,$ 
then the derivative of $f$ along the unit circle is continuous.
\end{thm}

\section{Modulus of semiannulus}

In the study of regularity for quasiconformal mappings (or more general
homeomorphisms),
the notion of ring domain (annulus) plays an important role.
In order to study the boundary regularity for homeomorphisms, 
we need its counterpart for the boundary.
In this section, we introduce the notion of semiannulus and present its
basic properties.
See also a survey article \cite{Sugawa10JA} by the third author
for expository accounts.

In what follows, we will consider subsets of the Riemann sphere
$\sphere=\C\cup\{\infty\}.$
Therefore, it is sometimes convenient to introduce the spherical (chordal)
distance
$$
\sd(z,w)=\frac{|z-w|}{\sqrt{(1+|z|^2)(1+|w|^2)}}.
$$
The spherical diameter of a set $A$ will be denoted by
$\sdiam A.$

A subset $S$ of $\sphere$ is called a {\it semiannulus}
if it is homeomorphic to
$$
T_R=\{z\in\C: 1\le |z|\le R, \Im z>0\}
$$
for some $R\in(1,+\infty).$
The two simple open arcs in the boundary of $S$ which
correspond to $\{|z|=1, \Im z>0\}$
and  $\{|z|=R, \Im z>0\}$ are called the {\it sides} of $S.$
The complementary boundary components of $S$ are called
the {\it ends} of $S.$

A semiannulus $S$ in a plane domain $D$
is said to be {\it properly embedded} in $D$
if $S\cap K$ is compact whenever $K$ is a compact subset of $D.$

A semiannulus $S$ is said to be {\it conformally equivalent} to
another semiannulus $S'$
if there is a homeomorphism $f:S\to S'$ which is conformal in $\Int S.$
We define the modulus of $S$ somewhat artificially as follows:
Set $\mod S=\log R$ when $S$ is conformally equivalent to $T_R.$
If there is no such an $R>1,$ we set $\mod S=0.$

One can also define $\mod S$ in an intrinsic way.
For that, we use the extremal length $\lambda(\Gamma)$
of a curve family $\Gamma$ (see \cite{Ahlfors:qc} for the definition
and its fundamental properties).

%Let $S$ be a semiannulus.
%properly embedded in the unit disk $\D$
Let $\Gamma_S$ be the collection of open arcs in $S$ dividing the two
sides of $S$ and $\Gamma_S'$ be that of closed arcs in $S$ joining the
two sides of $S.$
Here, a curve $\gamma$ in $S$ is called dividing if the sides of $S$
are contained in different connected components of $S\setminus\gamma.$
Then we have the following.

\begin{lem}\label{lem:char}
Let $S$ be a semiannulus in $\sphere.$
Then 
$$
\mod S=\frac{\pi}{\lambda(\Gamma_S)}=\pi\lambda(\Gamma_S').
$$
Furthermore, $\mod S=0$ if and only if there exists a sequence of 
simple closed arcs
$\gamma_n,~n=1,2,3,\dots,$ in $S$ joining the two sides of $S$
such that $\sdiam \gamma_n\to 0$ as $n\to\infty.$
\end{lem}

\begin{pf}
Let $h:\Int S\to\D$ be a Riemann mapping function.
By Carath\'eodory's theory of prime ends (see \cite{Pom:bound} for details),
the prime ends of $\Int S$ correspond to the boundary points of $\D$
through the function $h$ in a one-to-one fashion.
The sides of $S$ are open simple arcs in $\partial S,$
and therefore, these correspond to
disjoint open circular arcs, say, $O_1$ and $O_2$ in $\partial\D$
and $h$ extends to a homeomorphism of $S$ onto 
the semiannulus $S'=\D\cup(O_1\cup O_2).$
Let $C_1$ and $C_2$ be the connected components of 
$\partial\D\setminus(O_1\cup O_2).$
Clearly, $h(\Gamma_S)=\Gamma_{S'}$ is the collection of open arcs in $S'$
dividing the sides $O_1$ and $O_2.$
On the other hand, $h(\Gamma_S')=\Gamma_{S'}'$ is the collection of closed arcs
in $S'$ joining $O_1$ and $O_2.$
When we regard $S'$ as a quadrilateral, the curve families $\Gamma_{S'}$
and $\Gamma_{S'}'$ are conjugate to each other and therefore
$\lambda(\Gamma_{S'})=1/\lambda(\Gamma_{S'}')$ (see \cite{Ahlfors:qc}).
Since the extremal length is conformally invariant, we have
$\lambda(\Gamma_S)=1/\lambda(\Gamma_S').$

Obviously, $\lambda(\Gamma_S)=\lambda(\Gamma_{S'})=+\infty$ if and only if
either $C_1$ or $C_2$ reduces to one point, equivalently,
one of the two ends of $S$ is a prime end.
By the definition of prime ends, this means existence of a sequence of
simple closed arcs $\gamma_n$ (called a null-chain) in $S$ with 
the following properties:
$\gamma_n$ joins the two sides of $S,$
$\gamma_n$ separates $\gamma_{n-1}$ from $\gamma_{n+1}$ in $S,$
and $\sdiam\gamma_n\to0.$

Finally, we prove $\mod S=\pi/\lambda(\Gamma_S).$
This is certainly true when $\mod S=0$ by the above observation.
Thus we can assume that $\mod S>0.$
Then there exists a number $0<a<+\infty$
such that $S'=\D\cup O_1\cup O_2$ is conformally equivalent to
the rectangle $S_0=\{x+iy: 0\le x\le a, 0<y<\pi\}.$
Furthermore, the function $e^z$ maps  $S_0$ onto $T_{e^a}$ conformally inside.
Therefore, we now have $\mod S=\log e^a=a$ by definition.
On the other hand,
as is well known (cf.~\cite{Ahlfors:qc}), we have 
$\lambda(\Gamma_{S_0})=\pi/a.$
Hence, $\mod S=a=\pi/\lambda(\Gamma_{S_0})=\pi/\lambda(\Gamma_S).$
\end{pf}

In particular, if the two sides of $S$ have a positive spherical distance,
$\mod S>0.$
If $S$ is properly embedded in a good enough domain, then this sort of
condition is indeed characterizing positivity of the modulus.
For instance, we can show the following.

\begin{cor}\label{cor:mod}
Let $S$ be a semiannulus properly embedded in the unit disk $\D$
and let $U_1$ and $U_2$ be the connected components of $\D\setminus S.$
Then $\mod S>0$ if and only if the Euclidean distance between
$U_1$ and $U_2$ is positive.
\end{cor}

\begin{pf}
First note that $\dist(U_1,U_2)=\dist(\sigma_1,\sigma_2),$ where
``dist" stands for the Euclidean distance and
$\sigma_1, \sigma_2$ are the sides of $S.$
Since the spherical distance is comparable with the Euclidean distance
in $\D,$ the conclusion follows from the above lemma.
\end{pf}

The following result will constitute the basis of boundary estimates
for a disk homeomorphism.

\begin{thm}\label{thm:hyp}
Let $S$ be a semiannulus properly embedded in $\D$
and $U_1$ and $U_2$ be the two connected components of $\D\setminus S.$
Then
$$
\min\{\diam U_1, \diam U_2\}\le C\exp(-\tfrac12\mod S),
$$
where $C=4e^{\pi/2}.$
\end{thm}

For the proof of the theorem, we have to prepare a couple of lemmas.
Let $S$ be a semiannulus properly embedded in the unit disk $\D.$
If $\mod S$ is positive, by Lemma \ref{lem:char}, the Euclidean distance
between the two sides of $S$ is positive.
Therefore, the interior of the set
$
\overline S\cup\{1/\bar z: z\in S\}
$
becomes a ring domain and will be denoted by $\hat S.$
The modulus of a ring domain $B$ is defined to be $\log R$
if $B$ is conformally equivalent to the standard annulus $1<|z|<R.$
We now have the following by a symmetry principle.

\begin{lem}\label{lem:ref}
Let $S$ be a semiannulus properly embedded in the unit disk.
Then $\mod S=\mod \hat S$ whenever $\mod S>0.$
\end{lem}

\begin{pf}
Let $R=\exp(\mod S)>1.$
By definition, there is a homeomorphism $f:S\to T_R$
which is conformal on $\Int S.$
By the Schwarz reflection principle, $f|_{\Int S}$ can be continued
analytically to a conformal homeomorphism
$f:\hat S\to \{w\in\C: 1< |w|< R\}.$
Therefore, $\mod \hat S=\log R=\mod S.$
\end{pf}

We recall a sort of separation lemma for an annulus, which can date
back to Teichm\"uller's work in 1930's.
The following sharp form is due to Avkhadiev and Wirths \cite{AW05p}
(see also \cite[Theorem 3.17]{AW:sp} or \cite{Sugawa10JA}).

\begin{lem}[Avkhadiev-Wirths]\label{lem:sub}
Let $B$ be a ring domain in $\C$ with $\mod B>\pi$ which separates 
a given point $z_0\in\C$ from $\infty.$
Then there is a ring domain $A$ contained in $B$ of the form
$\{z: r<|z-z_0|<R\}$ such that $\mod A=\log\frac Rr\ge\mod B-\pi.$
The constant $\pi$ is sharp.
\end{lem}

A subset $S_0$ of a semiannulus $S$ is called a {\it subsemiannulus} of $S$
if $S_0$ is a semiannulus satisfying $\Gamma_{S_0}\subset\Gamma_S.$
Since $\lambda(\Gamma_{S_0})\ge\lambda(\Gamma_S),$
we have $\mod S_0\le\mod S$ by Lemma \ref{lem:char}.

For $\zeta\in\partial\D$ and $0<r_1<r_2<+\infty,$ we set
\be\label{eq:T}
T(\zeta; r_1,r_2)=\{z\in\D: r_1\le |\tfrac{z-\zeta}{z+\zeta}|\le r_2\}.
\ee
%A semiannulus in $\D$ of this form will be called {\it canonical}.
Note that $\hat T(\zeta; r_1,r_2)
=\{z\in\sphere: r_1< |\tfrac{z-\zeta}{z+\zeta}|< r_2\}$ and
$$
\mod T(\zeta; r_1,r_2)
=\mod\hat T(\zeta; r_1,r_2)=\log\frac{r_2}{r_1}.
$$

The following result is a hyperbolic analog of Lemma 2.7 in \cite{GMSV05}.

\begin{lem}\label{lem:hyp}
Let $T$ be a semiannulus properly embedded in $\D$ 
whose sides are circular arcs perpendicular to $\partial\D$
and let $V_1$ and $V_2$ be the connected components of $\D\setminus T.$
Then
$$
\min\{\diam V_1,\diam V_2\}\le\frac2{\cosh(\frac12\mod T)}.
$$
Equality holds if and only if $T$ is of the form $T(\zeta; r,1/r)$
for some $\zeta\in\partial\D$ and $0<r<1.$
\end{lem}

\begin{pf}
We denote by $d_\Omega$ the hyperbolic distance on a hyperbolic domain
$\Omega.$

Let $C_1$ and $C_2$ be the sides of $T$ and
let $\delta$ denote the hyperbolic distance between $C_1$ and $C_2$
in $\D.$
There is a unique hyperbolic line $C$ in $\D$ such that the hyperbolic length of
$C\cap T$ is $\delta,$ in other words, $C\cap T$ is the hyperbolic
geodesic joining $C_1$ and $C_2.$

Let $\zeta_1$ and $\zeta_2$ be the endpoints of $C$
with $\zeta_j\in \partial V_j\cap\partial\D,~j=1,2.$
We now define a conformal homeomorphism $L$ of $\D$ onto
the right half-plane $H$ by 
$$
L(z)=\frac{\zeta_2+z}{\zeta_2-z}-\frac{\zeta_2+\zeta_1}{\zeta_2-\zeta_1}.
$$
Then $L(C)$ is the half-line $(0,+\infty)$ and
$L(C_1)$ and $L(C_2)$ are concentric circular arcs centered at the origin.
We denote by $r_1$ and $r_2$ the radii of those circles.
Then the hyperbolic distance between $V_1$ and $V_2$ in $\D$ 
can be computed by
$$
\delta=d_\D(V_1,V_2)=d_H(L(V_1),L(V_2))=\int_{r_1}^{r_2}\frac{dx}{2x}
=\frac12\log\frac{r_2}{r_1}=\frac12\mod T.
$$
Thus the problem now reduces to finding a configuration of two hyperbolic
half-planes $V_1$ and $V_2$
with a fixed hyperbolic distance such that the minimum
of their Euclidean diameters is maximal (namely, the worst case).
Such a configuration is attained obviously by the situation that
$V_2=-V_1.$
In this case, $C$ becomes a line segment passing through the origin.
By a suitable rotation, we may assume that $\zeta_1=1, \zeta_2=-1.$
Let $a>0$ be the number determined by $V_1\cap\R=(a,1).$
Since $0$ is the midpoint of the geodesic $[-a,a]$ joining
$V_1$ and $V_2,$ we have $\delta/2=d_\D(0,a)=\arctanh a$
and $a=\tanh(\delta/2).$
The disk automorphism (hyperbolic isometry)
$g(z)=(z+a)/(1+az)$ maps the hyperbolic half-plane
$\{z\in\D: \Re z>0\}$ onto $V_1.$
Therefore, we see that $g(i)$ and $g(-i)$ are the endpoints of
$C_1=\partial V_1\cap\D$ and thus
$\diam V_1=|g(i)-g(-i)|=2(1-a^2)/(1+a^2).$
Finally, we get the estimate in this case
$$
\diam V_j=2\,\frac{1-\tanh(\delta/2)^2}{1+\tanh(\delta/2)^2}
=\frac2{\cosh\delta}.
$$
Since $\delta=\frac12\mod T,$ the estimate is now shown.
The equality case is obvious from the above argument.
\end{pf}

We are now ready to prove Theorem \ref{thm:hyp}.

\begin{pf}[Proof of Theorem \ref{thm:hyp}]
When $\mod S\le \pi,$ the assertion trivially holds.
We now suppose that $\mod S>\pi.$
Take points $\zeta_j\in\overline{U_j}\cap\partial\D~(j=1,2)$
and let $L(z)=(z+\zeta_2)/(z-\zeta_2)-(\zeta_1+\zeta_2)/(\zeta_1-\zeta_2).$
Then, by Lemma \ref{lem:sub}, $L(\hat S)$ contains a ring domain
$A$ of the form $\{w: r_1< |w|< r_2\}$ with $\mod A=\log\frac{r_2}{r_1}
\ge \mod \hat S-\pi=\mod S-\pi.$
Set $T=L\inv(\overline A)\cap\D$ and
let $V_1, V_2$ be the two components of $\D\setminus T$ so that
$U_j\subset V_j~(j=1,2).$
By Lemma \ref{lem:hyp}, we have
\begin{align*}
\min\{\diam U_1, \diam U_2\}
&\le \min\{\diam V_1,\diam V_2\} \\
&\le \frac2{\cosh(\tfrac12\mod T)} \\
&<4\exp(-\tfrac12\mod T) \\
&<4\exp(-\tfrac12\mod S+\pi/2).
\end{align*}
\end{pf}

We also have a variant of Theorem \ref{thm:hyp} for the case of half-planes.
It is an important point that we do not lose one half of the modulus
in the estimate with the expense of a condition for the modulus.
Note that a prototype can be found at \cite[Lemma 2.8]{GS01} and 
\cite[Lemma 2.8]{GMSV05}.

\begin{thm}\label{thm:uhp}
Let $S$ be a semiannulus properly embedded in $\uhp$
and $U_1$ and $U_2$ be the two connected components of $\uhp\setminus S.$
If $\mod S>\pi$ and if
$U_2$ is unbounded, then for any point $t_0\in\overline{U_1}\cap\R,$
the inequality
$$
\sup_{z\in U_1}|z-t_0|\le C_1\dist(t_0, U_2)\exp(-\!\mod S)
$$
holds, where $C_1=e^{\pi}.$
\end{thm}

\begin{pf}
%We have nothing to show when $\mod S\le \pi.$
%Thus we assume that $\mod S>\pi.$
Let $\hat S$ be the ring domain obtained by reflecting $\Int S$ in $\R.$
Then $\mod\hat S=\mod S>\pi$ by the same argument as in Lemma \ref{lem:ref}.
Since $\hat S$ separates $t_0$ from $\infty,$ Lemma \ref{lem:sub} guarantees
existence of numbers $0<r<R<+\infty$ such that
$A=A(t_0; r,R)\subset\hat S$ and $\mod A\ge \mod S-\pi.$
Since $\dist(t_0,U_2)\ge R,$ we now have
$$
\sup_{z\in U_1}|z-t_0|
\le r=R\exp(-\!\mod A(t_0; r,R))
\le \dist(t_0, U_2)\exp(\pi-\!\mod S).
$$
\end{pf}

\section{Boundary continuity of a homeomorphism}

The following simple example shows that a self-homeomorphism
of the unit disk $\D$ (sometimes called a disk homeomorphism)
does not necessarily have a continuous extension to the boundary:
$$
f(re^{i\theta})=r\exp i(\theta-\log(1-r)).
$$
Here is a criterion of continuous extendibility of 
a disk homeomorphism to a boundary point.
We recall that $T(\zeta; r, R)$ is defined in \eqref{eq:T}.

\begin{prop}\label{prop:ext}
Let $f:\D\to\D$ be a homeomorphism and let $\zeta\in\partial\D.$
The mapping $f$ extends continuously to $\zeta$ if
$$
\lim_{r\to0+}\mod f(T(\zeta; r,R))=+\infty
$$
for some $R>0.$
\end{prop}

\begin{pf}
Let $U_r$ be the connected component of 
$\D\setminus T(\zeta; r,R)$ with $\zeta\in\overline{U_r}$ for $0<r<R$
and let $V_R$ be the other one, which does not depend on $r.$
Then the family of the sets $U_r,~0<r<R,$ constitutes a fundamental system
of neighbourhoods of $\zeta.$
Theorem \ref{thm:hyp} now yields
$$
\min\{\diam f(U_r), \diam f(V_R)\}
\le C\exp(-\tfrac12\mod f(T(\zeta; r,R))).
$$
By assumption, the last term tends to 0 as $r\to0+.$
Since $\diam f(V_R)$ is a fixed number, this implies that
$\diam \overline{f(U_r)}\to0$ as $r\to0.$
Therefore, the intersection $\bigcap_{0<r<R}\overline{f(U_r)}$
consists of a single point.
We can now assign this point as the extended value of $f$ at $\zeta$
so that $f$ has a continuous extension to $\zeta.$
\end{pf}

We remark that the converse is not true in the last proposition.
Indeed, consider the homeomorphism $f:\D\to\D$ determined by
$f(\bar z)=\overline{f(z)},~z\in\D$ and
$$
f(re^{i\theta})=r\exp i\pi(\theta/\pi)^{-\log(1-r)},
\quad 0\le \theta\le\pi,~0<r<1.
$$
Then, by construction, $f$ extends to $1$ continuously by setting
$f(1)=1.$
However, since $f(re^{i\theta})\to1$ as $r\to1-$ for any fixed $\theta$
with $|\theta|<\pi,$ the converse of the proposition does not hold
(see the proof of the next theorem).

If the assumption of the last proposition is true for every point
of a non-degenerate subinterval of $\partial\D,$ then the converse holds.
As an application of our previous observations, we show indeed the
following theorem.

\begin{thm}\label{thm:gen}
Let $E$ be a subset of $\partial\D$ and
$f:\D\to\D$ be a homeomorphism.
Suppose that for every $\zeta\in E,$
$$
\lim_{r\to0+}\mod f(T(\zeta; r,R))=+\infty
$$
holds true for some number $R=R(\zeta)>0.$
Then $f$ extends to a continuous injective mapping of $\D\cup E$
into $\overline\D.$
\end{thm}

Since $\overline\D$ is a compact Hausdorff space, the inverse mapping
of a continuous bijection of $\overline\D$ onto itself is also continuous.
Therefore, as an immediate corollary, 
we have the following result of Brakalova \cite{Bra07}.
Note that, earlier than it, 
Jixiu Chen, Zhiguo Chen and Chengqi He \cite{CCH96}
proved a similar result in a special situation (see also the proof of
Lemma 2.3 in \cite{Chen01}).

\begin{cor}[Brakalova \cite{Bra07}]\label{cor:bra}
A homeomorphism $f:\D\to\D$ admits a homeomorphic extension to
$\overline\D$ if and only if for each $\zeta\in\partial\D,$ 
there is an $R=R(\zeta)>0$ such that
$$
\lim_{r\to0+}\mod f(T(\zeta; r,R))=+\infty.
$$
\end{cor}

\begin{pf}[Proof of Theorem \ref{thm:gen}]
By Proposition \ref{prop:ext}, $f$ can be extended continuously
to every point in $E.$ 
It is true even in the context of General Topology that
the extended mapping $f:\D\cup E\to\overline\D$ is indeed
continuous.

We next show that $f$ is injective in $E.$
First we observe that $f(E)\subset\partial\D.$
Therefore, if we assume that the conclusion does not hold, then
$f(\zeta_1)$ and $f(\zeta_2)$
are the same point, say, $\omega_0,$ for some
$\zeta_1,\zeta_2\in E$ with $\zeta_1\ne\zeta_2.$
We may take $R_j=R(\zeta_j)$ so small that 
$T(\zeta_1; r_1, R_1)\cap T(\zeta_2; r_2, R_2)=\emptyset$
for $0<r_1<R_1, 0<r_2<R_2.$
Let $U_1$ and $U_2$ be the connected components of
$\D\setminus T$ with $\zeta_1\in \overline{U_1},$
where $T=T(\zeta_1; r_1, R_1)$

Take sequences $z_n, z_n'\in\D,~n=1,2,\dots,$ 
so that $z_n\to \zeta_1$ and $z_n'\to \zeta_2$ as $n\to\infty.$
Since $z_n\in U_1$ and $z_n'\in U_2$ for a sufficiently large $n,$ one has
$\dist(f(U_1), f(U_2))\le |f(z_n)-f(z_n')|.$
We now let $n\to\infty$ to obtain $\dist(f(U_1), f(U_2))=0,$
which implies $\mod f(T)=0$ by Corollary \ref{cor:mod}.
This contradicts the condition $\mod f(T(\zeta_1; r_1,R_1))\to+\infty$
as $r_1\to0+.$
\end{pf}

In the same way, we have a half-plane version of the last theorem.

\begin{thm}\label{thm:uhpext}
A homeomorphism $f$ of the upper half-plane $\uhp$ admits
a homeomorphic extension to $\overline\uhp$ if and only if
for each $t\in\partial\uhp=\R\cup\{\infty\},$
$$
\lim_{r\to0+}\mod f(A(t;r,R)\cap\uhp)=+\infty
$$
for some $R=R(t)>0.$
\end{thm}

We recall that $A(\infty; r,R)$ is defined as $A(0; 1/R,1/r).$
We end the present section with an illustrating example.

\begin{example}\label{ex:1}
Firstly, consider the function $\tanh(\frac\pi4z),$
which maps the parallel strip $0<\Im z<1$ conformally onto 
$\{w\in\D: \Im w>0\},$ the upper half of the unit disk.
Observe that the line $\Im z=1$ is mapped by the function to the semicircle
$|w|=1, \Im w>0.$
We now define a self-homeomorphism $f$ of the upper half-plane
$\uhp$ by
$$
f(x+iy)=
\begin{cases}
\tanh\tfrac\pi4(x+iy) & \quad(0<y\le 1), \\
y \tanh\tfrac\pi4(x+i) & \quad(1<y).
\end{cases}
$$
Then $f:\uhp\to\uhp$ is a homeomorphism and conformal
in $0<\Im z<1.$
By construction, $f$ is not continuous at $\infty$ but satisfies
the assumptions of Theorem A except for the
continuity at $\infty.$
This example shows also that we cannot replace $\R\cup\{\infty\}$
by $\R$ in Theorem \ref{thm:uhpext}.
\end{example}

\section{Proof of main theorems}

In this section, we will give proofs of our main theorems stated in Section 1.
As we will see soon, we can show even slightly stronger assertions.
We need to introduce some technical quantities.

For a Beltrami coefficient $\mu,$
$$
K_\mu=\frac{1+|\mu|}{1-|\mu|}
$$
is sometimes called the {\it pointwise maximal dilatation} of $\mu$
(or, of a $\mu$-conformal homeomorphism $f$).
The $\mu$-conformal homeomorphism $f$ is $K$-quasiconformal precisely when
$K_\mu\le K$ a.e.
The pointwise maximal dilatation
is useful to measure quasiconformality of a $\mu$-conformal homeomorphism.
However, it is occasionally necessary to look at not only the absolute value
but also the argument of $\mu.$
For a Beltrami coefficient $\mu$ on a domain $\Omega$ and
a point $z_0\in\C$ (not necessarily in the domain $\Omega$), set
$$
D_{\mu,z_0}(z)=
\frac{\left|1-\mu(z)\frac{\bar z-\bar z_0}{z-z_0}\right|^2}{ 1-|\mu(z)|^2}
=\frac{|1-e^{-2i\theta}\mu(z)|^2}{1-|\mu(z)|^2}
$$
for $z\in\Omega,$ where $\theta=\arg(z-z_0).$
This quantity is sometimes called the directional dilatation
and it was introduced by Andreian Cazacu \cite{Cazacu57}.
This notion was effectively used by Reich and Walczak \cite{RW65}, 
Lehto \cite{Lehto68} and later by Brakalova and Jenkins \cite{BJ94},
\cite{BJ02}, Brakalova \cite{Bra10a}, \cite{Bra10b},
the first and third authors \cite{GS01}, 
Martio and the first author \cite{GM03}, and
Martio, Vuorinen and the authors \cite{GMSV05}.

It is easy to verify the inequalities
$$
\frac1{K_\mu(z)}\le D_{\mu,z_0}(z)\le K_\mu(z),\quad z\in\Omega.
$$
For a Beltrami coefficient $\mu$ on the upper half-plane $\uhp,$ we consider
the quantity
$$
Q_\mu(t; r,R)=\frac{1}{\pi\log(R/r)}
\iint_{A(t;r,R)\cap\uhp}\frac{D_{\mu,t}(z)}{|z-t|^2}dxdy
$$
for $t\in\R,$ and $0<r<R<+\infty.$

In terms of $Q_\mu(t; r,R)$ and $D_{\mu,t},$ 
we have distortion estimates for the modulus 
of a semiannulus under the $\mu$-conformal homeomorphism.
The following are variants of Proposition 2.4 and
Corollary 2.13 in \cite{GS01}
(see also \cite[Lemma 2.5]{GM03} for \eqref{eq:D}).
We omit the proof because we can show them in the same way as
in \cite{GS01}.

\begin{lem}
Let $\mu$ be a Beltrami coefficient on the upper half-plane $\uhp$
and $f$  be a $\mu$-conformal homeomorphism of $\uhp$ onto
another domain.
Then, for the semiannulus $T=A(t; r,R)\cap\uhp$
\be\label{eq:Q}
\frac1{Q_\mu(t; r,R)}\le\frac{\mod f(T)}{\mod T}\le Q_{-\mu}(t; r,R)
\ee
and
\be\label{eq:D}
-\frac1\pi\iint_T\frac{D_{-\mu,t}(z)-1}{|z-t|^2}dxdy
\le \mod T-\mod f(T)
\le \frac1\pi\iint_T\frac{D_{\mu,t}(z)-1}{|z-t|^2}dxdy.
\ee
\end{lem}

By making use of the last lemma, we are now able to prove 
Lemma \ref{lem:infty}.

\begin{pf}[Proof of Lemma \ref{lem:infty}]
We first note that
$$
Q_\mu(0; r_0, R)-1=\frac{2}{\pi\log(R/r_0)}
\iint_{A(0; r_0,R)\cap\uhp}
\frac{|\mu(z)|^2-\Re \frac{\bar z}{z}\mu(z)}{1-|\mu(z)|^2}
\frac{dxdy}{|z|^2}.
$$
Hence, the condition in the lemma means that
$Q_\mu(0; r_0, R)=1+o(\log R)=o(\log R)$ as $R\to+\infty.$
By \eqref{eq:Q}, we have
$$
\mod f(A(0; r_0,R)\cap\uhp)\ge \frac{\log(R/r_0)}{Q_\mu(0; r_0,R)}.
$$
The last term blows up when $R\to+\infty.$
Thus we now apply Theorem \ref{thm:uhpext} to obtain the desired conclusion.
\end{pf}

We are ready to prove Theorem \ref{thm:1}.

\begin{pf}[Proof of Theorem \ref{thm:1}]
First we observe the relations
$$
\frac{D_{\mu,t}(z)+D_{-\mu,t}(z)}2-1=\frac{2|\mu(z)|^2}{1-|\mu(z)|^2}
$$
and
$$
\frac{D_{\mu,t}(z)-D_{-\mu,t}(z)}2
%=-\frac{2e^{-2i\theta}\mu(z)}{1-|\mu(z)|^2}
=\frac{-2|z-t|^2}{1-|\mu(z)|^2}\Re\frac{\mu(z)}{(z-t)^2}.
$$
Therefore, conditions (1) and (2) in Theorem \ref{thm:1} mean
that the functions
\be\label{eq:int}
G(r)=\iint_{A(t;r,R_0)\cap\uhp}\frac{D_{\mu,t}(z)-1}{|z-t|^2}dxdy
\aand
H(r)=\iint_{A(t;r,R_0)\cap\uhp}\frac{D_{-\mu,t}(z)-1}{|z-t|^2}dxdy
\ee
both have finite limits as $r\to0+.$
(We remark here that similar conditions appear in \cite[Theorem 1.3]{Bra10b}
in connection with conformality condition at a point.)
By \eqref{eq:D}, we obtain 
$\mod f(A(t; r,R_0)\cap\uhp)=\log\frac1r+O(1)$ when $r\to0+.$
In particular, $\mod f(A(t;r,R_0)\cap \uhp)\to+\infty$ as $r\to0+.$
Therefore, by Proposition \ref{prop:ext}, we see that $f$ can be extended
to a homeomorphism of $\uhp\cup I$ into $\overline\uhp.$

From now on, we assume that $f(I)\subset\R.$
Again, by \eqref{eq:D}, we have
$$
H(r_1)-H(r_2)
\le \mod \big(f(A(t; r_1,r_2)\cap\uhp\big)-\log\frac{r_2}{r_1}
\le G(r_2)-G(r_1).
$$
Since $G(r)$ and $H(r)$ have finite limits as $r\to0+,$
\be\label{eq:asym1}
\lim_{r_2\to0+}\sup_{r_1\in(0,r_2)}
\left(\mod \big(f(A(t; r_1,r_2)\cap\uhp\big)-\log\frac{r_2}{r_1}\right)=0.
\ee
Therefore, by applying the argument of the proof of Lemma 4.1 in \cite{BJ94}
to a nesting sequence of semiannuli, we can show that the non-zero finite limit
$$
g(t):=\lim_{z\to t\text{~in~}\uhp}\left|\frac{f(z)-f(t)}{z-t}\right|
$$
exists for each $t\in I.$
Since the limit is unrestricted as long as $\Im z>0,$
we can let $\Im z\to0$ to obtain
$$
g(t)=\lim_{x\to t\text{~in~}\R}\left|\frac{f(x)-f(t)}{x-t}\right|
=\lim_{x\to t\text{~in~}\R}\frac{f(x)-f(t)}{x-t},
$$
which is nothing but the derivative of $f(t).$

When the convergence in (1) and (2) of the theorem is locally uniform for 
$t\in I,$ convergence in \eqref{eq:asym1} is also locally uniform for 
$t\in I.$
We can show now continuity of $f'(t)$ by the argument same as in \cite{BJ02}.
\end{pf}

The above proof was done along the same line as in \cite{BJ94} and \cite{BJ02}.
Therefore, we omitted details if the same argument in those papers
can be applied.
We should, however, note that the regularity assumptions on $\mu$-conformal
homeomorphisms are not clearly stated in their papers.
At least, their proof can be justified when $f$ belongs to the Sobolev
space $W^{1,1}_\loc,$
which is equivalent to that $f$ is an ACL homeomorphism
with locally integrable partial derivatives (see 
\cite[\S 4.9.2]{EG:mt} or \cite[Theorem 2.1.4]{Ziemer:wdf} for instance).
There is still a point where we should be careful.
We may extend the above $f:\uhp\cup I\to\overline\uhp$ 
furthermore to the lower half-plane $\uhp^-=\{z:\Im z<0\}$
by setting $f(z)=\overline{f(\bar z)}.$
Unlike the case of quasiconformal mappings, it is not clear that
the extended $f$ belongs to the Sobolev space 
$W^{1,1}_\loc(\uhp\cup I\cup\uhp^-).$
We give a simple example below to show that this is not necessarily
true in a general setting.
Therefore, the results in \cite{BJ94} might not be applied directly
under the current situations.
However, as we suggested in the proof, the arguments in \cite{BJ94} work
when we replace annuli by semiannuli in a suitable way.

\begin{example}
Define a function $\phi:\R\to\R$ by $\phi(y)=y\sin(1/y)$ if $0<|y|<1/\pi$
and $\phi(y)=0$ otherwise.
Then the function
$$
f(x+iy)=x+\phi(y)+iy, \quad x, y\in\R,
$$
is a self-homeomorphism of $\C$ satisfying $\overline{f(z)}=f(\bar z).$
Since the restriction of $f$ to $\uhp$ is a self-diffeomorphism,
$f\in W^{1,1}_\loc(\uhp).$
However, the image of the segment $[x, x+i/\pi]$ is of infinite length
for every $x\in\R,$
and therefore, $f$ is not ACL.
In particular, $f$ does not belong to $W^{1,1}_\loc(\C).$
\end{example}

\begin{rem}
In order to obtain convergence of the argument of $(f(z)-f(t))/(z-t),$
we shall need the convergence of the integrals as in \eqref{eq:int}
over the sets $\{z: r\le |z-t|\le R_0, \theta_1<\arg(z-t)<\theta_2\}$
for $0<r<R_0, 0\le \theta_1<\theta_2\le\pi$
(see \cite[\S 5]{BJ94}).
This condition is implied by the assumptions of Theorem A
but not by those of Theorem \ref{thm:1}.
\end{rem}

\begin{pf}[Proof of Theorem \ref{thm:2}]
Take a fixed point $z_0\in\uhp$ with $\Im z_0>R$ and set $w_0=f(z_0).$
We first note that condition \eqref{eq:t2} can be expressed by
$$
\frac1\pi\iint_{A(t;r,R)\cap\uhp}\frac{D_{\mu,t}(z)-1}{|z-t|^2}dxdy
\le \frac{2M}{\pi}
$$
for $t\in I$ and $0<r<R.$
Therefore, by \eqref{eq:D},
$$
\mod f(A(t;r,R)\cap\uhp)\ge \log \frac Rr-\frac{2M}\pi.
$$
In particular, $\mod f(A(t;r,R)\cap\uhp)>\pi$ for $0<r<r_0,$
where $r_0$ is taken so that $\log R/r_0\ge\pi+2M/\pi.$
Since $f(A(t; r, R)\cap\uhp)$ separates $w_0$ from $f(t)$ in $\uhp,$ one has
$\dist(f(t), V_2)\le |f(t)-w_0|,$ where $V_2$ is the unbounded component
of $\uhp\setminus f(A(t; r,R)\cap\uhp).$
We now take an arbitrary point $z\in\uhp$ with $|z-t|<r_0$
and set $r=|z-t|.$
Theorem \ref{thm:uhp} now yields
$$
|f(z)-f(t)|\le C_1|f(t)-w_0|\exp(-\log\tfrac Rr+\tfrac{2M}\pi)
=C_2|f(t)-w_0||z-t|,
$$
where $C_2=C_1e^{2M/\pi}/R.$
Thus we have shown that $f$ is locally Lipschitz continuous on $I.$
\end{pf}

\begin{pf}[Proof of Theorem \ref{thm:3}]
We set
$$
\omega(t; r)=\frac2{\pi r^2}\iint_{A(t;0, r)\cap\uhp}D_{\mu,t}(z)dxdy-1
=\frac2{\pi r^2}\iint_{A(t;0, r)\cap\uhp}(D_{\mu,t}(z)-1)dxdy.
$$
Then \eqref{eq:t3} is equivalent to the condition
$\limsup_{r\to0+}\omega(t; r)\le \alpha\inv-1.$
Since the convergence is locally uniform in $I$ by assumption,
for a compact subset $I_0$ of $I$
and a given $0<\alpha'<\alpha,$ we can find 
a constant $R>0$ such that $\omega(t; r)\le 1/\alpha''-1$ for
$t\in I_0$ and $0<r\le R,$ where $\alpha''=(\alpha+\alpha')/2.$
On the other hand, we observe that $\omega(t; r)\ge -1$ by definition.
In particular, $\omega(t; r)$ is bounded in $0<r\le R$ for each $t\in I_0.$
We now have the relation
\begin{align*}
(Q_\mu(t; r,R)-1)\log \frac Rr
&=\frac{1}{\pi}\iint_{A(t;r,R)\cap\uhp}\frac{D_{\mu,t}(z)-1}{|z-t|^2}dxdy \\
&=\frac{\omega(t; R)-\omega(t; r)}2+\int_r^R\omega(t; s)\frac{ds}s
\end{align*}
(see \cite[Lemma 3.8]{GS01}).
In particular,
$$
(Q_\mu(t; r, R)-1)\log\frac Rr
\le \frac1{2\alpha''}
+\int_r^R\left(\frac1{\alpha''}-1\right)\frac{ds}{s}
=\frac1{2\alpha''}+\left(\frac1{\alpha''}-1\right)\log\frac Rr,
$$
and consequently, $Q_\mu(t; r,R)\le 1/\alpha'$ for $0<r<r_0$ and $t\in I_0$
for a sufficiently small $0<r_0<R.$
By \eqref{eq:Q}, we have
$$
\mod f(A(t; r,R)\cap\uhp)
\ge \frac{\mod (A(t; r,R)\cap\uhp)}{Q_\mu(t; r,R)}
\ge \alpha'\log\frac Rr
$$
for $0<r<r_0.$
As in the proof of Theorem \ref{thm:2}, we have an estimate of the form
$|f(z)-f(t)|\le C|z-t|^{\alpha'}$ for $z\in\uhp$ $t\in I_0$ with $|z-t|\le r_0.$
\end{pf}

\begin{pf}[Proof of Theorem \ref{thm:1d}]
For a fixed $\zeta=e^{i\theta_0}\in\partial\D,$ we define a M\"obius transformation $L$
by $L(z)=i(\zeta-z)/(\zeta+z), \tilde L(z)=i(f(\zeta)-z)/(f(\zeta)+z)$ and let $M=L\inv.$
Note that $L$ and $\tilde L$ map $\D$ conformally onto $\uhp$ and that $L(\zeta)=0,
\tilde L(f(\zeta))=0.$
Let $\hat\mu$ be a Beltrami coefficient on $\uhp$ given by
$\hat\mu=\mu\circ M\cdot \overline{M'}/M'.$ 
Then $F=\tilde L\circ f\circ M$ is a $\hat\mu$-conformal self-homeomorphism of $\uhp.$
For $w=u+iv=L(z),$ we therefore have the relation
$$
\frac{\hat\mu(w)dudv}{(w-0)^2(1-|\hat\mu(w)|^2)}
=\frac{L'(z)^2}{L(z)^2}\frac{\mu(z)dxdy}{(1-|\mu(z)|^2)}.
$$
Since $L'(z)/L(z)=2\zeta/(z^2-\zeta^2),$ we see that conditions (i) and (ii)
in Theorem \ref{thm:1d} are equivalent to conditions (1) and (2) in Theorem \ref{thm:1}
with $t=0,$ respectively.
The continuity of $f$ now follows from the argument same as in the proof of Theorem
\ref{thm:1} and from Proposition \ref{prop:ext}.
If we write $e^{i\theta}=M(u)$ for $u\in\R,$ we have
$ie^{i\theta}d\theta=M(u)du.$
In particular, $d\theta/du=2$ at $u=0.$
Since $F=\tilde L\circ f\circ M$ and $\tilde L'(0)=1/(2if(\zeta)),$ we have
$$
\lim_{u\to0}\frac{F(u)}{u}=\frac1{if(\zeta)}\lim_{\theta\to\theta_0}
\frac{f(e^{i\theta})-f(e^{i\theta_0})}{\theta-\theta_0},
$$
which implies the differentiability property of $f.$
The continuity of $df(e^{i\theta})/d\theta$ follows from that of Theorem \ref{thm:1}.
\end{pf}

\def\cprime{$'$} \def\cprime{$'$} \def\cprime{$'$}
\providecommand{\bysame}{\leavevmode\hbox to3em{\hrulefill}\thinspace}
\providecommand{\MR}{\relax\ifhmode\unskip\space\fi MR }
\providecommand{\MRhref}[2]{%
  \href{http://www.ams.org/mathscinet-getitem?mr=#1}{#2}
}
\providecommand{\href}[2]{#2}

%\bibliography{papers}
\end{document}